\documentclass[12pt]{article}

\usepackage{amsmath,amssymb,amsbsy,amsfonts,amsthm,latexsym,
            amsopn,amstext,amsxtra,euscript,amscd}

\newtheorem{lem}{Lemma}[section]
\newtheorem{lemma}[lem]{Lemma}
\newtheorem{remark}[lem]{Remark}

\newtheorem{proposition}[lem]{Proposition}

\newtheorem{theorem}[lem]{Theorem}

\newtheorem{corollary}[lem]{Corollary}

\def\\{\cr}
\def\({\left(}
\def\){\right)}
\def\[{\left[}
\def\]{\right]}
\def\<{\langle}
\def\>{\rangle}

\def\N{{\mathbb N}}

\def\eps{\varepsilon}

\begin{document}

\title{On the maximal order of numbers in the ``factorisatio numerorum'' 
problem}

\author{
{Martin Klazar} \\
{\normalsize  Department of Applied Mathematics (KAM)} \\
{\normalsize  and Institute for Theoretical Computer Science (ITI), Charles 
University} \\
{\normalsize Malostransk\'e n\'am\v est\'\i\ 25, 118 00 Praha, 
Czech Republic}\\
{\normalsize {\tt klazar@kam.mff.cuni.cz}} \\
\and
{Florian~Luca} \\
{\normalsize Instituto de Matem{\'a}ticas, Universidad Nacional
Aut\'onoma de M{\'e}xico} \\
{\normalsize Ap. Postal 61-3 (Xangari), C.P.~58089, Morelia, Michoac{\'a}n, 
M{\'e}xico} \\
{\normalsize {\tt fluca@matmor.unam.mx}}}

\date{\today}

\pagenumbering{arabic}

\date{\today}
\maketitle

\begin{abstract}
Let $m(n)$ be the number of ordered factorizations of $n\ge 1$ in factors larger 
than $1$. We prove that for every $\eps>0$
$$
m(n)<\frac{n^{\rho}}{\exp\((\log n)^{1/\rho}/(\log\log n)^{1+\eps}\)}
$$
holds for all integers $n>n_0$, while, for a constant $c>0$, 
$$
m(n)>\frac{n^{\rho}}{\exp\(c(\log n/\log\log n)^{1/\rho}\)}
$$
holds for infinitely many positive integers $n$, where $\rho=1.72864\dots$ is 
the real solution to $\zeta(\rho)=2$. We investigate also arithmetic properties of
$m(n)$ and the number of distinct values of $m(n)$.
\end{abstract}

\section{Introduction}

Let $m(n)$ be the number of ordered factorizations of a positive integer $n$ 
in factors bigger than $1$. For example, $m(12)=8$ since we have the factorizations 
$12$, $2\cdot 6$, $6\cdot 2$, $3\cdot 4$, $4\cdot 3$, $2\cdot 2\cdot 3$, 
$2\cdot 3\cdot 2$, and $3\cdot 2\cdot 2$. By the definition, $m(1)=0$ but we 
will see that in some situations it is useful to set $m(1)=1$ or $m(1)=1/2$.
Kalm\'{a}r \cite{kalm31a} found the average order of $m(n)$: for
$x\to\infty$, 
\begin{equation}
\label{eq:F}
M(x)=\sum_{n\le x} m(n)=cx^{\rho}(1+o(1)),
\end{equation}
where $\rho=1.72864\dots$ is the real solution to $\zeta(\rho)=2$  and 
$c=0.31817\dots$ is given by $c=-1/\rho\zeta'(\rho)$. (As usual, $\zeta(s)=
\sum_{n\ge 1}n^{-s}$.) Further detailed and strong results 
on the average order of $m(n)$ were obtained by Hwang \cite{hwan00}.

\medskip

In contrast, good bounds on the maximal order of $m(n)$ were lacking.
Erd\H os claimed in the end of his article \cite{erdo41} that there exist 
positive constants $0<c_2<c_1<1$ such that 
$$
m(n)<\frac{n^{\rho}}{\exp\((\log n)^{c_2}\)}
$$
holds for all $n>n_0$, while
$$
m(n)>\frac{n^{\rho}}{\exp\((\log n)^{c_1}\)}
$$
holds for infinitely many $n$, but he gave no details. To our knowledge, the best proved
bounds on the maximal order state that $m(n)<n^{\rho}$ for every $n\ge 1$ (Chor, Lemke 
and Mador \cite{chor00}, a simple proof by induction was recently given by Coppersmith and Lewenstein 
\cite{copp_lewe}) and that for any $\eps>0$ one has $m(n)>n^{\rho-\eps}$ for infinitely many $n$ 
(Hille \cite{hill37}, \cite{copp_lewe} gives an explicit construction). (In Lemma~\ref{lem:up} we 
strengthen the argument of \cite{chor00} and show that $m(n)\le n^{\rho}/\sqrt{2}$ for every $n\ge 1$.)

\medskip

Here we come close to determining the maximal order of $m(n)$.  We prove that it is, roughly, 
$n^{\rho}/\exp\((\log n)^{1/\rho}\)$. More precisely, we prove that for every $\eps>0$, 
$$
m(n)<\frac{n^{\rho}}{\exp\((\log n)^{1/\rho}/(\log\log n)^{1+\eps}\)}
$$
holds for all $n>n_0$ (Theorem~\ref{upbound}), while 
$$
m(n)>\frac{n^{\rho}}{\exp\(c(\log n)^{1/\rho}/(\log\log n)^{1/\rho}\)}
$$
holds with certain constant $c>0$ for infinitely many positive integers $n$ (Theorem~\ref{lobound}). 

\medskip

The paper is organized 
as follows. In Section 2 we give auxiliary results, of which Lemma~\ref{lem:riemann} on the speed of convergence 
$\rho_k\to\rho$ ($\rho_k$ is a ``finite'' counterpart of $\rho$ for $m(n)$ restricted to smooth numbers $n$ with 
no prime factor exceeding $p_k$, the $k$th prime number) and Lemmas~\ref{lem:up}--\ref{lemma:down} 
giving explicit inequalities for $m(n)$ and $m_k(n)$ ($m_k(n)=m(n)$ if $n$ has no prime factor $>p_k$ and 
$m_k(n)=0$ else) may be of independent interest. Section 3 is devoted to the proof of the upper bound. The proof 
is elementary (uses real analysis only) and is obtained by combining the combinatorial bounds on $m(n)$ 
in Lemmas~\ref{lem:up} and \ref{lem:comb}, standard bounds from the prime numbers theory, and 
the convergence bound in Lemma~\ref{lem:riemann}. Section 4 is devoted to the proof of the lower bound. In the first
version of this article, still available at \cite[version 1]{KL}, we proved by an elementary 
approach similar to that in Section 3, with the additional ingredience being Kalm\'ar's asymptotic relation 
(\ref{eq:F}), a weaker lower bound that has $(\log n)^{1/\rho}$ in the denominator replaced with 
the bigger power $(\log n)^{\rho/(\rho^2-1)+o(1)}$. Here, we prove in Section 4 a lower bound with the matching
exponent $1/\rho$ of the $\log n$ by a method 
suggested to us by an anonymous referee. The method works in the complex domain and combines the uniform (i.e., with 
error estimates independent on $k$) version of (\ref{eq:F}) for $m_k(n)$, bounds on smooth numbers, 
and again Lemma~\ref{lem:riemann}. In Section 5, we give further references and comments on the 
history of $m(n)$ and some related problems. We also investigate arithmetical properties of $m(n)$ and prove,
 for example, that $m(n)$ is not eventually periodic modulo $k$ for any integer 
$k>1$, and that $m(n)$ is not a holonomic sequence.

\bigskip

\noindent{\bf Acknowledgments.}
Most of this paper was written during a very enjoyable visit
by the first author to the Mathematical Institute 
of the UNAM in Morelia, Mexico, in March 2005. This author wishes to express his thanks to 
that institution for the hospitality and support. He also acknowledges the 
support to ITI by the project 1M0021620808 of the Czech	Ministry of Education. Both authors 
are deeply grateful to an anonymous referee who outlined for them a plan how to prove a better
lower bound, which is now carried out in Section 4. 

\section{Preliminaries and auxiliary results}

Let us begin with recalling some notation. For a positive integer $n$ we write $\omega(n)$ and 
$\Omega(n)$ for the number of distinct prime factors of $n$ and the total number of 
prime factors of $n$ (including multiplicities), respectively. 
We use the letters $p$ and $q$ with or without subscripts to denote prime numbers.
We put $P(n)$ for the largest prime factor of $n$. We write $\log $ for the 
natural logarithm. In the complex domain 
(mainly in Section 4) we use $s$ to denote generic variable and write $\sigma$ and $\tau$ for its real and imaginary
part, respectively, so $s=\sigma+i\tau$, where $i=\sqrt{-1}$. We use 
the Vinogradov symbols $\ll$ and $\gg$ and the Landau symbols $O$ and $o$ 
with their usual meanings.

\medskip

The proof of the following estimate is standard and we omit it.

\begin{lemma}
\label{lem:primeestimate}
If $\delta>\delta_0>1$, then the estimate
\begin{equation}
\label{eq:deltalow}
\sum_{p>t} \frac{1}{p^{\delta}} = \frac{(\delta-1)^{-1}}{t^{\delta-1}\log t}+
O\left(\frac{1}{t^{\delta-1}(\log t)^2}\right)
\end{equation}
holds uniformly for $t>2$.
\end{lemma}

\medskip
Let $p_k$ be the $k$th prime. We shall use the well known asymptotic relations
$$
\sum_{p\le x}\log p=x+O(x/\log x)
$$
(equivalent to the Prime Number Theorem) and
$$
p_k=k\log k+k\log\log k+O(k)
$$
(the full asymptotic expansion $p_k=k(\log k+\log\log k-1+\cdots)$ was found by Cipolla \cite{cipo}).
Let ${\cal P}_k$ be the set (including $1$) of 
positive integers composed only of the primes $p_1=2,p_2,\dots,p_k$, and 
$m_k(n)$ be the number of ordered factorizations of $n$ in factors lying in 
${\cal P}_k\backslash\{1\}$. We allow $k=\infty$, then $p_k=\infty$, 
${\cal P}_{\infty}={\cal P}$ is the set of all primes, and 
$m_{\infty}(n)=m(n)$. Note that, for $k\in\N$,  $m_k(n)>0$ iff $n\in{\cal P}_k$, if 
$m_k(n)>0$ then $m_k(n)=m(n)$, and if $n\le p_k$ then $m_k(n)=m(n)$. Let, for complex $s$ with $\sigma>1$ 
and $k\in\N\cup\{\infty\}$,
$$
\zeta_k(s)=\prod_{p\le p_k}\left(1-\frac{1}{p^{s}}\right)^{-1}=
\sum_{n\in{\cal P}_k}\frac{1}{n^s}
$$
and $\rho_k$ be the real solution to $\zeta_k(\rho_k)=2$. For $k=\infty$ we get
the Euler-Riemann zeta function $\zeta(s)=\zeta_{\infty}(s)$ and the number $\rho=\rho_{\infty}$.
Note that for $k\in\N$ the series for $\zeta_k(s)$ converges absolutely even for $\sigma>0$. 
For every $s$ with $\sigma>1$ we have the convergence $\zeta_k(s)\to\zeta(s)$ as $k\to\infty$. 
For $k\in\N\cup\{\infty\}$, one has the identity (setting $m_k(1)=1$ for every $k$)
$$
\sum_{n\ge 1}\frac{m_k(n)}{n^s}=\sum_{l\ge 0}(\zeta_k(s)-1)^l=\frac{1}{2-
\zeta_k(s)},
$$
which implies that $m_k(n)=o(n^{\rho_k+\eps})$ for every fixed $\eps>0$. Our approach to estimating $m(n)$ is based on
approximating the ``infinite'' quantities $m(n)$, $\rho$, and $\zeta(s)$ with their ``finite'' 
counterparts $m_k(n)$, $\rho_k$, and $\zeta_k(s)$ for $k\in\N$ but $k\to\infty$. We quantify the degrees of 
approximation in the following two lemmas. The first lemma is obtained by considering the infinite series
defining $\zeta_k(s)$ and $\zeta(s)$ and its easy proof is omitted.

\begin{lemma}\label{zeta_k_and_zeta}
We have 
$$
\rho_1=1<\rho_2=1.43527\dots<\rho_3=1.56603\dots<\dots<\rho=1.72864\dots
$$ 
and $\rho_k\to\rho$ as $k\to\infty$. 
The convergence $\zeta_k(s)\to\zeta(s)$ as $k\to\infty$ is uniform on every complex domain 
$\sigma>\sigma_0>1$ and the same is true for the convergence $\zeta_k'(s)\to\zeta'(s)$ and 
for all higher derivatives. Also, for every $k\in\N\cup\{\infty\}$ we have $\zeta_k'(\rho_k)<0$.
\end{lemma}

\noindent
We shall use this lemma to bound various expressions containing $\rho_k$, $\zeta_k(\rho_k)$, 
$\zeta_k(s)$, $1/\zeta_k'(\rho_k)$ etc. by constants independent on $k$.

\begin{lemma}
\label{lem:riemann}
The estimate 
$$
\rho-\rho_k=\frac{1}{(\rho-1)|\zeta'(\rho)|}\cdot 
\frac{1}{k^{\rho-1}(\log k)^{\rho}}\left(1+O\left(\frac{\log\log k}{\log k}
\right)\right)
$$
holds for all $k\ge 2$.
\end{lemma}

\begin{proof} We will assume that $k\ge 2$.
The equation $\zeta_k(\rho_k)^{-1}=\zeta(\rho)^{-1}=1/2$ implies that
$$
\prod_{2\le p\le p_k}\left(1-\frac{1}{p^{\rho_k}}\right)=\prod_{p\ge 2} 
\left(1-\frac{1}{p^{\rho}}\right).
$$
Taking logarithms and regrouping, we get
\begin{equation*}
\sum_{2\le p\le p_k}\left(\log\(1-\frac{1}{p^{\rho}}\)- 
\log\(1-\frac{1}{p^{\rho_k}}\)\right)=-\sum_{p>p_k}
\log\(1-\frac{1}{p^{\rho}}\).
\end{equation*}
The left side satisfies, by Lagrange's Mean-Value Theorem (the derivative of the function 
$x\mapsto \log(1-1/p^x)$ is $(\log p)/(p^x-1)$), 
\begin{eqnarray}
\label{eq:intvaltheorem}
\sum_{2\le p\le p_k}\log\(1-\frac{1}{p^{\rho}}\)- 
\log\(1-\frac{1}{p^{\rho_k}}\)&=&(\rho-\rho_k)\sum_{2\le p\le p_k}\frac{\log p}{p^{\sigma_p}-1}\\
&>&(\rho-\rho_k)(\log 2)/3\nonumber
\end{eqnarray}
for some numbers $\sigma_p\in(\rho_k,\rho)\subset (1.4,1.8)$. The right side is 
\begin{eqnarray}
\label{eq:right}
-\sum_{p>p_k}\log\(1-\frac{1}{p^{\rho}}\) & = & 
\sum_{p>p_k}\frac{1}{p^{\rho}}+O\(\sum_{p>p_k} \frac{1}{p^{2\rho}}\)\nonumber\\
& = & \frac{(\rho-1)^{-1}}{p_k^{\rho-1}\log(p_k)}
\left(1+O\left(\frac{1}{\log k}\right)\right)\nonumber\\
& = & \frac{(\rho-1)^{-1}}{k^{\rho-1}(\log k)^{\rho}}\left(1+O\left(
\frac{\log\log k}{\log k}\right)\right),
\end{eqnarray}
where we used Lemma~\ref{lem:primeestimate} and the fact that $p_k=k(\log k+O(\log\log k))$. We
get immediately that
\begin{equation}
\label{eq:rho-rhok}
\rho-\rho_k\ll\frac{1}{k^{\rho-1}(\log k)^{\rho}}.
\end{equation}

\medskip
To do better, we return to (\ref{eq:intvaltheorem}) and write 
$$
\frac{\log p}{p^{\sigma_p}-1}  
=  
\frac{\log p}{p^{\rho}-1}\left(1+\frac{p^{\sigma_p}}{p^{\sigma_p}-1}
(p^{\rho-\sigma_p}-1)\right).
$$
We have $1\le p^{\sigma_p}/(p^{\sigma_p}-1)\le 2$ and, using (\ref{eq:rho-rhok}),
$$
p^{\rho-\sigma_p}-1\le\exp((\rho-\rho_k)\log p_k)-1\ll (\rho-\rho_k)\log p_k
\ll\frac{1}{k^{\rho-1}(\log k)^{\rho-1}}.
$$
Hence the right side of (\ref{eq:intvaltheorem}) equals
\begin{eqnarray*}
(\rho-\rho_k)\sum_{2\le p\le p_k}\frac{\log p}{p^{\sigma_p}-1}
& = & (\rho-\rho_k)(1+O(k^{1-\rho}(\log k)^{1-\rho}))\sum_{2\le p\le p_k}\frac{\log p}{p^{\rho}-1}\\
& = & (\rho-\rho_k)(1+O(k^{-1/2}))\sum_{2\le p\le p_k}\frac{\log p}{p^{\rho}-1}. 
\end{eqnarray*}
Equating the right side of (\ref{eq:intvaltheorem}) and (\ref{eq:right}) we get the relation
$$
(\rho-\rho_k)\sum_{2\le p\le p_k}\frac{\log p}{p^{\rho}-1}=
\frac{(\rho-1)^{-1}}{k^{\rho-1}(\log k)^{\rho}}\left(1+O\left(
\frac{\log\log k}{\log k}\right)\right).
$$
All is left to notice is that 
\begin{eqnarray*}
|\zeta'(\rho)| & = & \sum_{p\ge 2}\frac{\log p}{p^{\rho}-1}=
\sum_{p\le p_k}\frac{\log p}{p^{\rho}-1}+\sum_{p>p_k}\frac{\log p}{p^{\rho}-1}
\\
& = & \sum_{p\le p_k}\frac{\log p}{p^{\rho}-1}+O(k^{-1/2}),
\end{eqnarray*}
where the last estimate follows again from Lemma~\ref{lem:primeestimate} via the fact that 
$\log p\ll p^{1/10}$:
$$
\sum_{p>p_k}\frac{\log p}{p^{\rho}-1}\ll \sum_{p>p_k}
\frac{1}{p^{\rho-0.1}}\ll \frac{1}{p_k^{\rho-1.1}\log p_k}<
1/\sqrt{k}.
$$
The claimed estimate now follows.
\end{proof}

In the next three lemmas, we prove combinatorial inequalities involving $m_k(n)$ and $m(n)$. In the first 
lemma, we slightly improve the result from \cite[Theorem 5]{chor00} that $m_k(n)<n^{\rho_k}$ for every $n\ge 1$.
The second lemma is crucial for obtaining bounds of the type $m(n)=o(n^{\rho})$. The third lemma gives some 
lower estimates on $m(n)$.

\begin{lemma}
\label{lem:up}
For every $k\in\N\cup\{\infty\}$ and  $n\ge 1$ (with $m_k(1)=0$),
$$
m_k(n)\le \frac{1}{\sqrt{2}}\;n^{\rho_k}.
$$
\end{lemma}

\begin{proof}
For every $r,s\ge 1$ we have (now setting $m_k(1)=0$)
\begin{equation}\label{ineqwith2}
m_k(rs)\ge 2m_k(r)m_k(s).
\end{equation}
To show this inequality, we assume that $r,s\ge 2$ (for $r=1$ or $s=1$ it 
holds trivially) and consider the set $X$ of all pairs $(u,v)$ where $u$ ($v$) 
is an ordered factorization of $r$ ($s$) in factors lying in 
${\cal P}_k\backslash\{1\}$, and 
the set $Y$ of the same factorizations of $rs$. If $u$ is 
$r=d_1\cdot d_2\cdot\ldots\cdot d_i$ and $v$ is $s=e_1\cdot e_2\cdot\ldots
\cdot e_j$, we define the factorizations of $rs$
\begin{eqnarray*}
F((u,v))&=&d_1\cdot d_2\cdot\ldots\cdot d_i\cdot e_1\cdot e_2\cdot\ldots\cdot 
e_j\\
G((u,v))&=&d_1\cdot d_2\cdot\ldots\cdot d_{i-1}\cdot (d_ie_1)\cdot e_2\cdot
\ldots\cdot e_j.
\end{eqnarray*}
The inequality (\ref{ineqwith2}) follows from the fact that the mappings $F$ 
and $G$ 
are injections from $X$ to $Y$ which moreover have disjoint images. We leave a 
simple verification of this fact to the reader.

Suppose now that $m_k(n_0)>n_0^{\rho_k}/\sqrt{2}$ for some $n_0\ge 2$. By 
(\ref{ineqwith2}),  
we have $m_k(n_0^2)\ge 2m_k(n_0)^2>n_0^{2\rho_k}$ and hence we can take some
$\eps>0$ so that $m_k(n_0^2)\ge (n_0^2)^{\rho_k+\eps}$. Then, again by (\ref{ineqwith2}), 
$m_k(n_0^{2i})\ge (n_0^{2i})^{\rho_k+\eps}$ for every $i=1,2,\ldots$, which is in 
contradiction with $m_k(n)=o(n^{\rho_k+\eps})$.
\end{proof}

\begin{lemma}
\label{lem:comb}
Suppose that $q_1,\dots,q_k$ are primes, not necessarily distinct, such that 
the product $q_1q_2\dots q_k$ divides $n$. Then, with $m(1)=1$, 
\begin{equation}
\label{eq:combfact}
m(n)< (2\Omega(n))^k\cdot m(n/q_1q_2\dots q_k).
\end{equation}
\end{lemma}

\begin{proof}
It suffices to prove only the case $k=1$; i.e., the inequality
\begin{equation}\label{clearineq}
m(n)< 2\Omega(n)\cdot m(n/p),
\end{equation}
where $p$ is a prime dividing 
$n$, because the general case follows easily by iteration. 
Let $X$ be the set of all pairs $(u,i)$ where $u$ is an ordered factorization 
of $n/p$ (in parts bigger than $1$) and $i$ is an integer satisfying $1\le i\le 2r+1$, where $r$ is the number 
of parts in $u$. Let $Y$ be the set 
of all ordered factorizations of $n$ in parts bigger than $1$. 
We shall define a surjection $F$ from 
$X$ onto $Y$. This will prove 
(\ref{clearineq}) because $r\le\Omega(n/p)=\Omega(n)-1$, and therefore for 
every $u$ we have 
$2r+1<2\Omega(n)$ pairs $(u,i)$, and so 
$$
m(n)=|Y|\le|X|<2\Omega(n)\cdot m(n/p). 
$$
For 
$(u,i)\in X$, where $u$ is 
$n/p=d_1\cdot d_2\cdot\ldots\cdot d_r$, we define $j=i-r$ and set $F((u,i))$ 
to be the 
factorization
$$
n=d_1\cdot\ldots\cdot d_{i-1}\cdot(pd_i)\cdot d_{i+1}\cdot \ldots\cdot d_r
$$ 
if $1\le i\le r$ and 
$$
n=d_1\cdot\ldots\cdot d_{j-1}\cdot p\cdot d_j\cdot\ldots\cdot d_r
$$ 
if $r+1\le i\le 2r+1$ (for $j=1$, $p$ is the first part and for $j=r+1$ it is 
the last one). 
It is clear that $F$ is a surjection.
\end{proof}

\begin{lemma}\label{lemma:down}
If $n_1,n_2,\dots,n_k$ are positive integers such that for no $i\ne j$ 
we have $n_i|n_j$, then
$$
m(n_1n_2\dots n_k)\ge k!\cdot m(n_1)m(n_2)\dots m(n_k).
$$
This implies that for every $n\ge 1$ we have
$$
m(n)\ge\omega(n)!\cdot 
2^{\Omega(n)-\omega(n)}\ \mbox{ and }\ m(n)\ge 2^{\Omega(n)-1}.
$$
\end{lemma}

\begin{proof}
Let $X$ be the set of all $k$-tuples $(u_1,u_2,\dots,u_k)$, where $u_i$ is an 
ordered factorization of 
$n_i$ in parts bigger than $1$ and let $Y$ be the set of these factorizations 
for $n_1n_2\dots n_k$. For every 
permutation $\sigma$ of $1,2,\dots,k$, we define a mapping 
$F_{\sigma}:\;X\to Y$ by
$$
F_{\sigma}((u_1,u_2,\dots,u_k))=u_{\sigma(1)}\cdot u_{\sigma(2)}
\cdot\ldots\cdot u_{\sigma(k)},
$$ 
i.e., we concatenate factorizations $u_i$ in the order prescribed by 
$\sigma$. It is clear that each 
$F_{\sigma}$ is an injection. Suppose that 
$F_{\sigma}((u_1,u_2,\dots,u_k))=F_{\tau}((v_1,v_2,\dots,v_k))$
for some permutations $\sigma,\tau$ and factorizations $u_i$ and $v_i$. 
It follows that $u_{\sigma(1)}$ is 
an initial segment of $v_{\tau(1)}$ or vice versa, and hence $n_{\sigma(1)}$ 
divides $n_{\tau(1)}$ or vice versa.
This implies that $\sigma(1)=\tau(1)$ and $u_{\sigma(1)}=v_{\tau(1)}$. 
Applying the same argument, we obtain that 
$\sigma(j)=\tau(j)$ and $u_{\sigma(j)}=v_{\tau(j)}$ also for $j=2,\dots,k$. 
Thus $\sigma=\tau$ and 
$u_j=v_j$ for $j=1,2,\dots,k$. We have proved that the $k!$ mappings 
$F_{\sigma}$ have mutually disjoint images. 
Therefore 
$$
k!\cdot m(n_1)m(n_2)\dots m(n_k)=k!|X|\le |Y|=m(n_1n_2\dots n_k).
$$

If $n=q_1^{a_1}q_2^{a_2}\dots q_k^{a_k}$ is the prime factorization of $n$, 
applying the first inequality to 
the $k$ numbers $n_i=q_i^{a_i}$ and using that $m(p^a)=2^{a-1}$, we obtain
$$
m(n)\ge k!\prod_{i=1}^k 2^{a_i-1}=k!\cdot 2^{\Omega(n)-k}, 
$$
which is the second inequality. Using that $k!/2^k\ge 1/2$ for every $k\ge 1$, 
we get the third inequality. 
\end{proof}

\noindent
Note that $m(n)\ge 2^{\Omega(n)-1}$ is tight for every $n=p^a$.

\section{The upper bound}

We prove the following upper bound on the maximal order of $m(n)$.

\begin{theorem}\label{upbound}
For every $\eps>0$ we have 
$$
m(n)<\frac{n^{\rho}}{\exp\((\log n)^{1/\rho}/(\log\log n)^{1+\eps+o(1)}\)}
$$
for integers $n>2$.
\end{theorem}
\begin{proof} Let $\eps>0$ be given. To bound $m(n)$ from above, 
we split the integers $n>0$ in two groups, those with $\omega(n)\le k$ and those with
$\omega(n)>k$, which we shall treat by different arguments; the optimum value of the parameter $k=k(n)$  
will be selected in the end.

{\em The case $\omega(n)\le k$.}
Let $n=q_1^{a_1}q_2^{a_2}\ldots q_r^{a_r}$, $r\le k$, be the prime 
decomposition of $n$ where $q_1<q_2<\ldots<q_r$. We denote by $\bar{n}$ 
the number obtained from $n$ by replacing $q_i$ in the 
decomposition by $p_i$, the $i$th smallest prime. Then $\bar{n}\le n$. From the fact that 
$m(n)$ depends only on the exponents $a_i$ and from Lemma~\ref{lem:up} we get 
$$
m(n)=m(\bar{n})=m_r(\bar{n})<\bar{n}^{\rho_r}\le n^{\rho_k}.
$$
Thus, by Lemma~\ref{lem:riemann}, 
\begin{eqnarray}
\label{eq:case1}
m(n) & < & n^{\rho_k}\nonumber\\
& = & n^{\rho}\exp(-(\rho-\rho_k)\log n)\nonumber\\
& = & n^{\rho}\exp\left(-(c+o(1))\frac{\log n}{k^{\rho-1}(\log k)^{\rho}}
\right)
\end{eqnarray}
where $c=(\rho-1)^{-1}|\zeta'(\rho)|^{-1}>0$.

{\em The case $\omega(n)>k$.} Let $l(n)$ be the product of some $k$ distinct prime factors of $n$; then 
$l(n)\ge p_1p_2\dots p_k$, the product of the $k$ smallest primes. We have the estimates
$$
\sum_{p\le p_k} \log p=p_k+O(p_k/\log p_k)=k\log k+k\log\log k+O(k)
$$
and
$$
2\Omega(n)\le (2/\log 2)\log n<3\log n.
$$
By Lemmas~\ref{lem:comb}, \ref{lem:up} and these estimates,
\begin{eqnarray}
\label{eq:case2}
m(n) & < & (2\Omega(n))^{k} m(n/\ell(n)) 
 <  (3\log n)^{k}\frac{n^{\rho}}{\ell(n)^{\rho}}\nonumber\\
& \le & (3\log n)^{k}\frac{n^{\rho}}{(p_1\dots p_{k})^{\rho}}\nonumber\\
& = & n^{\rho}\exp\big(-k(\rho\log k+\rho\log\log k-\log\log n+O(1))\big).
\end{eqnarray}

To determine the best upper bound on $m(n)$, we begin with $k$ in the form $k=k(n)=(\log n)^{\alpha+o(1)}$ 
where $\alpha\in(0,1)$ is a constant. Necessarily $\alpha\ge 1/\rho$, for else the argument of 
$\exp$ in (\ref{eq:case2}) 
is eventually positive and we get a useless bound. It follows that the optimum is 
$\alpha=1/\rho$ when the arguments of both $\exp$s in 
(\ref{eq:case1}) and (\ref{eq:case2}) are $-(\log n)^{1/\rho+o(1)}$, provided that 
\begin{equation}\label{positivity}
\rho\log k+\rho\log\log k-\log\log n+O(1)>c>0 
\end{equation}
for big $n$. Now we set, more precisely, 
$$
k=k(n)=\frac{(\log n)^{1/\rho}}{(\log\log n)^{d+o(1)}}
$$
with a constant $d>0$. With this $k$, the function in (\ref{positivity}) becomes 
$\rho(1-d+o(1))\log\log\log n+O(1)$ and we see that condition (\ref{positivity}) is satisfied for 
$d<1$ (for $d>1$ the argument of the $\exp$ in (\ref{eq:case2}) is again eventually positive). With
this $k$, the arguments of the $\exp$s in (\ref{eq:case1}) and (\ref{eq:case2}) are, respectively, 
$$
-\frac{(\log n)^{1/\rho}}{(\log\log n)^{1+(\rho-1)(1-d)+o(1)}}\ \mbox{ and }\ 
-\frac{(\log n)^{1/\rho}}{(\log\log n)^{d+o(1)}}.
$$
Setting $d=1-\eps/2(\rho-1)$, we obtain the stated bound. 
\end{proof} 

\section{The lower bound}

We prove the following lower bound on the maximal order of $m(n)$.
\begin{theorem}\label{lobound}
There is a constant $c>0$ such that the inequality
$$
m(n)>\frac{n^{\rho}}{\exp\(c(\log n/\log\log n)^{1/\rho}\)}
$$
holds for infinitely many integers $n>0$.
\end{theorem} 

\noindent
We shall see that it is possible to take $c=2.01630\dots-\eps$.
We begin with explaining the effective Ikehara--Ingham theorem on Dirichlet series. Then we apply it 
to $1/(2-\zeta_k(s))$ to obtain an asymptotic relation for the average order of $m_k(n)$ with error estimate
independent on $k$. Finally, combining this relation with an estimate on density of 
smooth numbers we obtain Theorem~\ref{lobound}. For the background on Dirichlet series we refer to Tenenbaum 
\cite{tene}.

\medskip

Suppose that $(a_n)_{n\ge 1}$ is a sequence of nonnegative real numbers with the summatory function
$$
A(t)=\sum_{n\le e^t}a_n
$$
and the Dirichlet series
$$
F(s)=\sum_{n=1}^{\infty}\frac{a_n}{n^s}=\int_{0-}^{\infty} e^{-st}\;\mathrm{d}A(t).
$$
Suppose that $F(s)$ converges for $\sigma>a>0$. We may assume that $a$ is the abscissa of (absolute) 
convergence; then by the Phragm\'en--Landau theorem, $a$ is a singularity of $F(s)$. The effective 
Ikehara-Ingham theorem, proved by Tenenbaum \cite{tene} (who used the method of Ganelius 
\cite{gane}), extracts an asymptotic relation for $A(x)$ as 
$x\to\infty$ from the local behavior of $F(s)$ near $a$ 
and, moreover, it provides an explicit estimate of the error term in terms of  
the regularity of $F(s)$ on the vertical segments $a+\sigma+i\tau$, $-T\le\tau\le T$, as $\sigma\to 0+$. 
We quote the theorem verbatim from Tenenbaum \cite[p. 234]{tene}.

\begin{theorem} {\bf (``Effective'' Ikehara-Ingham).}\label{iithm}
Let $A(t)$ be a non-decreasing function such that the integral 
$$
F(s):=\int_0^{\infty} e^{-st}\;\mathrm{d}A(t)
$$
converges for $\sigma>a>0$. Suppose that there exist constants $c\ge 0$,
$\omega>-1$, such that the function
$$
G(s):=\frac{F(s+a)}{s+a}-\frac{c}{s^{\omega+1}}\ \ (\sigma>0)
$$
satisfies
\begin{equation}\label{regcondition}
\eta(\sigma,T):=\sigma^{\omega}\int_{-T}^T|G(2\sigma+i\tau)-G(\sigma+i\tau)|
\;\mathrm{d}\tau=o(1)\ \ (\sigma\to 0+)
\end{equation}
for each fixed $T>0$. Then we have
\begin{equation}\label{asymrel}
A(x)=\left\{\frac{c}{\Gamma(\omega+1)}+O(\rho(x))\right\}e^{ax}x^{\omega}\ \ (x\ge 1),
\end{equation}
with
$$
\rho(x):=\inf_{T\ge 32(a+1)}\{T^{-1}+\eta(1/x,T)+(Tx)^{-\omega-1}\}.
$$
Furthermore, the implicit constant in (\ref{asymrel}) depends only on $a$, $c$, and
$\omega$. An admissible choice for this constant is 
$$
52+1652c(a+1)(\omega+1)+69c(1+(\omega+1)e^{1-\omega}(\omega+1)^{\omega+2})/\Gamma(\omega+1).
$$
\end{theorem}

\noindent
Note that for a meromorphic $F(s)$ with a simple pole at $s=a$ (so $\omega=0$), 
the condition (\ref{regcondition}) is satisfied iff $F(s)$ has on the line $\sigma=a$ no other poles. 

\medskip

We shall apply Theorem~\ref{iithm} to the functions
$$
F(s)=F_k(s)=\sum_{n\ge 1}\frac{m_k(n)}{n^s}=\frac{1}{2-\zeta_k(s)}
$$
for $k\ge 2$, $a=\rho_k$, $c=c_k=-1/\rho_k\zeta_k'(\rho_k)$, and $\omega=0$. It is not hard to 
prove (we do this in the next Proposition) that $\rho_k$ is the only pole of $F_k(s)$ on 
$\sigma=\rho_k$ when $k\ge 2$ (this is not true for $k=1$) and thus by Theorem~\ref{iithm} 
$$
\sum_{n\le x}m_k(n)=(c_k+o(1))x^{\rho_k}\ \ (x\to\infty)
$$
for each fixed $k\ge 2$. (In contrast, $\sum_{n\le x}m_1(n)=2^r-1$ where $2^r\le x<2^{r+1}$.)
To get a good lower bound on $m(n)$, we have to strengthen this by obtaining
uniformity in $k$ of the error term $o(1)$. This follows from Theorem~\ref{iithm}, once we prove that 
for $F(s)=F_k(s)$ the condition (\ref{regcondition}) is satisfied uniformly in $k$. 

\begin{proposition}\label{unifsatisf}
Let, for $k\ge 2$,
$$
G_k(s)=\frac{F_k(s+\rho_k)}{s+\rho_k}-\frac{c_k}{s}=
\frac{1}{(2-\zeta_k(s+\rho_k))(s+\rho_k)}-\frac{c_k}{s}
$$
and $T>0$ be arbitrary but fixed. Then
$$
\lim_{\sigma\to 0+}\int_{-T}^T|G_k(2\sigma+i\tau)-G_k(\sigma+i\tau)|
\;\mathrm{d}\tau=0
$$
uniformly in $k\ge 2$; that is, the condition (\ref{regcondition}) holds uniformly in $k$.
\end{proposition}

\begin{proof}
Let $t(\sigma)=\sigma^{1/5}$; any function $t(\sigma)>0$ satisfying, as $\sigma\to 0+$, that
$t(\sigma)\to 0$ and $\sigma/t(\sigma)^4\to 0$ would do in our argument. For every fixed $T>0$, we bound 
the integrand by a quantity that depends only on $\sigma$ and not on $\tau$ and $k\ge 2$ and that goes to $0$ 
as $\sigma\to 0+$; this will prove the statement. We manage doing this by splitting $[-T,T]$ in two ranges,
$t(\sigma)\le|\tau|\le T$ and $|\tau|\le t(\sigma)$, in which we apply different
arguments.

{\em The range $t(\sigma)\le |\tau|\le T$.} Denoting by $\gamma$ the horizontal segment with endpoints 
$\sigma+i\tau$ and $2\sigma+i\tau$, we have the bound 
$$
|G_k(2\sigma+i\tau)-G_k(\sigma+i\tau)|=\left|\int_{\gamma}G'_k(z)\;dz\right|\le
\sigma|G_k'(s_0)|
$$
where $s_0$ is some point lying on $\gamma$. The derivative of $G_k(s)$ equals
$$
G_k'(s)=
\frac{(s+\rho_k)\zeta_k'(s+\rho_k)
+\zeta_k(s+\rho_k)-2}{(2-\zeta_k(s+\rho_k))^2(s+\rho_k)^2}+\frac{c_k}{s^2}.
$$
We bound the numerators and denominators of this expression. As for the numerators, by 
Lemma~\ref{zeta_k_and_zeta}, there is a 
constant $c=c(T)>0$ depending only on $T$ such that 
$$
|(s+\rho_k)\zeta_k'(s+\rho_k)
+\zeta_k(s+\rho_k)-2|,\ |c_k|<c
$$
holds for every $k\ge 2$ and $s$ with $0<\sigma<1$ and $|\tau|\le T$. For the second denominator, 
we have, in our range and for $0<\sigma<1$,
$$
\frac{\sigma}{|s_0|^2}\le \frac{\sigma}{\sigma^2+t(\sigma)^2}=
\frac{\sigma^{3/5}}{\sigma^{8/5}+1}<\sigma^{3/5}.
$$
We bound the first denominator. Clearly, $|s+\rho_k|^2\ge \rho_k^2>1$ for every $s$ with $\sigma>0$. 
For every $k\ge 2$ and every $s$ with $0<\sigma<1$ and any $\tau$ we have
$$
|2-\zeta_k(s+\rho_k)|\ge {\text{\rm Re}}(2-\zeta_k(s+\rho_k))
=\sum_{\substack{n\ge 1\\ P(n)\le p_k}}\frac{1}{n^{\rho_k+\sigma}}
\left(n^{\sigma}-\cos(\tau\log n)\right)
$$
and, consequently, (recall that $k\ge 2$ and $1<\rho_k<2$)
$$
|2-\zeta_k(s+\rho_k)|^2>\left(\frac{2-\cos(\tau\log 2)-\cos(\tau\log 3)}{27}\right)^2=:h(\tau).
$$
Since $2^{\alpha}=3$ holds for no fraction $\alpha$, $h(\tau)=0$ only for $\tau=0$. The function 
$h(\tau)$ is continuous and even and $h(\tau)\sim\beta\tau^4$ as $\tau\to 0$ for a constant $\beta>0$. 
Thus there is a constants $\beta_1=\beta_1(T)<1$ depending only on $T$ such that if 
$0<\sigma<\beta_1$ then the minimum of $h(\tau)$ on $[t(\sigma),T]$ is attained at $t(\sigma)$ and 
$h(t(\sigma))>\beta t(\sigma)^4/2$. Hence, in our range and for $0<2\sigma<\beta_1$, 
$$
\frac{\sigma}{|2-\zeta_k(s_0+\rho_k)|^2\cdot|s_0+\rho_k|^2}<
\frac{2\sigma}{\beta t(\sigma)^4}=\frac{2\sigma^{1/5}}{\beta}.
$$
Taking together all estimates, we have in our range and for $0<\sigma<\beta_1/2$ that 
$$
|G_k(2\sigma+i\tau)-G_k(\sigma+i\tau)|\le\sigma|G_k'(s_0)|<c(2\sigma^{1/5}/\beta+\sigma^{3/5}),
$$
which is the required bound.

{\em The range $|\tau|\le t(\sigma)$. }We prove that there is an absolute constant $\delta>0$ such that
for every $k\ge 2$ and $s$ with $|s|<\delta$ we have the expansion
$$
G_k(s)=d_k+O(s),
$$
where $d_k$ is a constant and the constant implicit in $O$ is absolute. (We need independence on $k$ both 
for the constant in $O(s)$ and for the domain of validity of the error estimate.)
Then if $0<\sigma<\delta^5/32$ and 
$|\tau|\le t(\sigma)$, both numbers $\sigma+i\tau$ and $2\sigma+i\tau$ satisfy $|s|<\delta$, and we have
the bound
$$
|G_k(2\sigma+i\tau)-G_k(\sigma+i\tau)|=O(|\sigma+i\tau|+|2\sigma+i\tau|)=O(\sigma^{1/5})
$$
with absolute constants in $O$s, which is the required bound.

We begin with the origin-centered closed disc $B=B(0,0.1)$; the point of the radius $0.1$ is only that 
$\rho_2-0.1>1$. We define functions $f_k(s)$ by 
$$
f_k(s)=\frac{\zeta_k(s+\rho_k)-2-s\zeta_k'(\rho_k)-s^2\zeta_k''(\rho_k)/2}{s^3}.
$$
Let $a_k$ be the maximum value taken by $|\zeta_k(s)|$   
on the circle $|s-\rho_k|=0.1$. By the maximum modulus principle ($f_k(s)$ 
is holomorphic on $B$), for every $s\in B$ we have
$$
|f_k(s)|\le 10^3(a_k+2+10^{-1}\zeta_k'(\rho_k)+10^{-2}\zeta_k''(\rho_k)/2).
$$
Thus, by Lemma~\ref{zeta_k_and_zeta}, there is an absolute constant $M>0$ such that
$$
|f_k(s)|<M
$$
holds for every $s\in B$ and every $k\ge 2$. We rewrite 
$\zeta_k(s+\rho_k)=2+s\zeta_k'(\rho_k)+s^2\zeta_k''(\rho_k)/2+s^3f_k(s)$ 
as
\begin{eqnarray*}
\frac{1}{(2-\zeta(s+\rho_k))(s+\rho_k)} & = & -
\frac{1}{s\rho_k\zeta_k'(\rho_k)}\times \frac{1}{1+
s/\rho_k}\\
& \times & \frac{1}{1+s\zeta_k''(\rho_k)/2\zeta_k'(\rho_k)
+s^2 f_k(s)/\zeta_k'(\rho_k)}\\
&=&-\frac{1}{s\rho_k\zeta_k'(\rho_k)}\times \frac{1}{1+
s/\rho_k}\times\frac{1}{1+sb_k+s^2h_k(s)}.
\end{eqnarray*}
It follows, by Lemma~\ref{zeta_k_and_zeta} and the bound $|f_k(s)|<M$ valid on $B$, that there is 
a $\delta$, $0<\delta<0.1$, such that $|s/\rho_k|<1/2$ and $|sb_k+s^2h_k(s)|<1/2$ whenever 
$|s|<\delta$ and $k\ge 2$.
Using the estimate $(1+s)^{-1}=1-s+O(s^2)$, valid for $|s|<1/2$, and Lemma~\ref{zeta_k_and_zeta} 
we obtain for $k\ge 2$ and 
$|s|<\delta$ the expansion
\begin{eqnarray*}
\frac{1}{(2-\zeta_k(s+\rho_k))(s+\rho_k)} & = & 
\frac{c_k}{s}\left(1-\frac{s}{\rho_k}
+O(s^2)\right)\left(1-s\frac{\zeta_k''(\rho_k)}{2\zeta_k'(\rho_k)}
+O(s^2)\right)\\
& = & \frac{c_k}{s}-c_k\left(\frac{1}{\rho_k}+
\frac{\zeta_k''(\rho_k)}{2\zeta_k'(\rho_k)}\right)+O(s),
\end{eqnarray*}
where $c_k=-1/\rho_k\zeta_k'(\rho_k)$ and the constants in $O$s are absolute. Now the required expansion  
$G_k(s)=d_k+O(s)$ (valid for $|s|<\delta$ and with an absolute constant in $O$) is immediate.
\end{proof}

\begin{corollary}
\label{thm:uniform}
There is a constant $\beta_2>2$ such that for every $x>\beta_2$ and every $k\ge 2$ we have 
$$
\sum_{\substack{n\le x\\ P(n)\le p_k}} m(n)=\sum_{n\le x} m_k(n)
>x^{\rho_k}/5.
$$
\end{corollary}
\begin{proof}
By Theorem~\ref{iithm} and Proposition~\ref{unifsatisf}, there is a function $e(x)>0$ such that 
$e(x)\to 0$ as $x\to\infty$, and for every $x\ge 1$ and every $k\ge 2$ we have
$$
\left|\sum_{n\le x} m_k(n)-c_k x^{\rho_k}\right|<e(x)x^{\rho_k}.
$$
The sequence of $c_k=-1/\rho_k\zeta_k'(\rho_k)$, $k=1,2,\dots$, monotonically decreases and converges to 
$c_{\infty}=-1/\rho\zeta'(\rho)>0.3$. Thus if $x$ is big enough so that $e(x)<0.1$, the sum 
$\sum_{n\le x} m_k(n)$ must be bigger than $0.2x^{\rho_k}$.
\end{proof}

\medskip
We proceed to the proof of Theorem~\ref{lobound}. We denote, as usual, 
$$
\Psi(x,y)=\#\{n\le x:P(n)\le y\}.
$$
By Corollary \ref{thm:uniform}, for every $k\ge 2$ and $x>\beta_2$
there exists an $n_0\le x$ such that 
$$
\Psi(x,p_k)m(n_0)>x^{\rho_k}/5=\frac{x^{\rho}}{5\exp(
(\rho-\rho_k)\log x)}.
$$
We select $k=k(x)$ so that it satisfies
$$
k=(\log x)^{\alpha+o(1)}
$$ 
as $x\to\infty$, for some absolute constant $\alpha\in(0,1)$ (we make our choice of $k$ more precise 
later). Then 
$$
p_{k}=(1+o(1))k\log k=(\log x)^{\alpha+o(1)}. 
$$
A theorem  due to de Bruijn, see Theorem 2 in Tenenbaum's
book \cite[p. 359]{tene}, shows that
$$
\log(\Psi(x,p_{k}))=(1+o(1))Z,
$$
where 
\begin{eqnarray*}
Z & = & \frac{\log x}{\log p_{k}}\log\left(1+\frac{p_{k}}{\log x}\right)+
\frac{p_{k}}{\log p_{k}}
\log\left(1+\frac{\log x}{p_{k}}\right)\\
& = & \frac{p_{k}}{\log p_{k}}(1+o(1))+\frac{p_{k}}{\log p_{k}}
\log\left(1+\frac{\log x}{p_{k}}\right)\\
& = & (1+o(1))k(\log\log x-\log k).
\end{eqnarray*}
By Lemma \ref{lem:riemann},  
$$
\rho-\rho_k=\frac{c_1+o(1)}{k^{\rho-1}(\log k)^{\rho}} 
$$
where $c_1=1/((\rho-1)|\zeta'(\rho)|)$.  
Substituting both estimates in the lower bound on $\Psi(x,p_k)m(n_0)$, we get (absorbing the $5$ in 
the denominator in the $o(1)$ terms)
$$
m(n_0)>\frac{x^{\rho}}{\exp\left(c_1(1+o(1))
\frac{\log x}{k^{\rho-1}(\log k)^{\rho}}+
(1+o(1))k(\log\log x-\log k)\right)}.
$$
This suggests to choose $k$ so that both terms in the argument of the exponential, 
$$
\frac{\log x}{k^{\rho-1}(\log k)^{\rho}}\qquad {\text{\rm and }}\qquad 
k(\log\log x-\log k),
$$
are of the same order of magnitude. This occurs when $\alpha=1/\rho$, 
more precisely when 
$$
k=\lfloor d(\log x)^{1/\rho}(\log\log x)^{-(\rho+1)/\rho}\rfloor
$$
with any constant $d>0$, because then
$$
\frac{\log x}{k^{\rho-1}(\log k)^{\rho}}=d^{1-\rho}\rho^{\rho}(1+o(1))
\left(\frac{\log x}{\log\log x}\right)^{1/\rho}
$$
and 
$$
k(\log\log x-\log k)=(1-\rho^{-1})d (1+o(1))
\left(\frac{\log x}{\log\log x}\right)^{1/\rho}.
$$
Thus, for this selection of $k$,  
$$
m(n_0)> 
\frac{x^{\rho}}{\exp\left((c+o(1))\left(\frac{\log x}{\log\log x}\right)^{1/\rho}
\right)}
$$
where $c>0$ is a constant depending only on the choice of $d$. The lower bound eventually 
increases monotonically to infinity, and we conclude that there exist infinitely many numbers 
$n_0$ satisfying 
$$
m(n_0)>
\frac{n_0^{\rho}}{\exp\left((c+o(1))\left(\frac{\log n_0}{\log\log n_0}\right)^{1/\rho}
\right)}. 
$$
The proof of Theorem~\ref{lobound} is complete.

\medskip
It is not difficult to find the optimal value of $d$; it yields the value 
$$
c=(\rho^{\rho+1} c_1)^{1/\rho}=\left(\frac{\rho^{\rho+1}}{(\rho-1) 
|\zeta'(\rho)|}\right)^{1/\rho}\approx 2.01630.
$$

\section{Historical remarks and arithmetical properties of $m(n)$}

We begin with a survey of some previous results on $m(n)$. We restrict our attention 
only to works dealing directly with this quantity. There are many other variants of 
factorization counting functions
(with restrictions on factors, counting unordered factorizations etc.) and for a survey
on these we refer the reader to Knopfmacher and Mays \cite{knop_mays}. 

\medskip

Kalm\'ar proved in \cite{kalm31b} that the error term $o(1)$ in (\ref{eq:F}) is 
$$
O(\exp(-\alpha\log\log x\cdot\log\log\log x)), \mbox{ with } \alpha<
\frac{1}{2(\rho-1)\log 2}\approx 1.97996.
$$ 
Ikehara devoted three papers to the estimates of 
$M(x)$. In \cite{ikeh37}, he gave weak bounds of the type 
$M(x)>x^{\rho-\varepsilon}$ on a sequence of $x$ tending to infinity, 
and $M(x)<x^{\rho+\varepsilon}$ for all large enough $x$. In the review of 
\cite{ikeh37}, 
Kalm\'ar pointed out a gap in the proof and sketched a correct
argument. In \cite{ikeh39}, Ikehara gave a proof of (\ref{eq:F}) with an 
error bound $O(\exp(q\log\log x))$ for some constant $q<0$, which is slightly 
weaker than 
Kalm\'ar's result. Finally, in \cite{ikeh41}, he succeeded to get a stronger 
error bound
$$
O(\exp(-\alpha(\log\log x)^{\gamma})), \mbox{ with } \alpha>0\mbox{ and }
\gamma<4/3.
$$
Hwang \cite{hwan00} obtained an improvement of Ikehara's last bound by 
replacing $4/3$ 
with $3/2$.

\medskip

Rieger proved in \cite{rieg61}, besides other results, that for all positive 
integers 
$k,l$ with $(k,l)=1$ one has 
$$
\sum_{n\le x,\,n\equiv l\,(k)}m(n)=\frac{1+o(1)}{\varphi(k)}M(x)=
\frac{-1}{\varphi(k)\rho\zeta'(\rho)}\cdot x^{\rho}(1+o(1)).
$$
Warlimont investigated in \cite{warl93} variants of $m(n)$ counting ordered 
factorizations with distinct parts and with coprime parts and estimated their 
summatory functions. Hille in \cite{hill37} proved that 
$m(n)=O(n^{\rho})$ and that $m(n)>n^{\rho-\varepsilon}$ for infinitely many 
$n$. We already mentioned in Section 1 the remark of Erd\H{o}s on $m(n)$ in 
\cite{erdo41} 
and we mentioned (and improved) the result of Chor, Lemke 
and Mador \cite{chor00} that $m(n)<n^{\rho}$ for all $n$. Other elementary 
and constructive proofs of the bounds 
$m(n)\le n^{\rho}$ and $\limsup_n m(n)/n^{\rho-\epsilon}=\infty$ were recently 
given by Coppersmith and Lewenstein \cite{copp_lewe}.

\medskip

We now turn to recurrences and explicit formulas. The recurrence $m(1)=1$ and
\begin{equation}\label{recdivi}
m(n)=\sum_{d|n,\,d<n}m(d)\ \mbox{ for $n>1$}
\end{equation}
is immediate from fixing the first part in a factorization. If we set 
$m^*(1)=1/2$ and $m^*(n)=m(n)$ for $n>1$, then
$2m^*(n)=\sum_{d|n}m(d)$ holds for all $n\ge 1$. By
M\"obius inversion, $m(n)=2\sum_{d|n}\mu(d)m^*(n/d)$ for all $n\ge 1$. For 
$n=q_1^{a_1}q_2^{a_2}\dots q_r^{a_r}>1$ this can be rewritten as the recurrence formula
\begin{equation}
\label{Moebinv}
m(n)=2\left(\sum_{i}m\bigg(\frac{n}{q_i}\bigg)-
\sum_{i<j}m\bigg(\frac{n}{q_iq_j}\bigg)+\cdots+(-1)^{r-1}m\bigg(\frac{n}
{q_1q_2\dots q_r}\bigg)\right),
\end{equation}
in which we must set $m(1)=1/2$. 
Formulas (\ref{recdivi}) and (\ref{Moebinv}) are from Hille's paper \cite{hill37}.
In fact, (\ref{Moebinv}) is stated there incorrectly with 
$m(1)=1$, as was pointed out by K\"uhnel \cite{kuhn41} and Sen \cite{sen41}.

\medskip

Clearly, $m(p^a)=2^{a-1}$ because ordered factorizations of $p^a$ in parts 
$>1$ are in bijection with (additive) compositions of $a$ in parts $>0$. If 
$p\ne q$ are primes and $a\ge b\ge 0$ are integers, we have the formula
$$
m(p^aq^b)=2^{a+b-1}\sum_{k=0}^b\binom{a}{k}\binom{b}{k}2^{-k}
$$
that was derived in \cite{chor00} and before by Sen \cite{sen41} and 
MacMahon \cite{macm93}. In 
particular,
\begin{equation}\label{mforpq}
m(p^a q)=(a+2)2^{a-1}\ \mbox{ and }\ m(p^a q^2)=(a^2+7a+8)2^{a-2}.
\end{equation}
In general, for $n=q_1^{a_1}q_2^{a_2}\dots q_r^{a_r}$, and
$a=a_1+a_2+\cdots+a_r$, MacMahon
\cite{macm93} derived the
formula
$$
m(q_1^{a_1}q_2^{a_2}\dots
q_r^{a_r})=\sum_{j=1}^a\sum_{i=0}^{j-1}(-1)^i\binom{j}{i}
\prod_{k=1}^r \binom{a_k+j-i-1}{a_k}.
$$
A more complicated summation formula for $m(q_1^{a_1}q_2^{a_2}\dots 
q_r^{a_r})$ but involving only nonnegative summands was obtained by 
K\"uhnel in \cite{kuhn41} and \cite{kuhn50}. Let $d_k(n)$ be the number of 
solutions of 
$n=n_1n_2\dots n_k$, where $n_i\ge 1$ are positive integers; so $d_2(n)$ is 
the number of 
divisors of $n$. Sklar \cite{skla52} mentions the formula 
\begin{equation}\label{Sklarova}
m(n)=\sum_{k=1}^{\infty}\frac{d_k(n)}{2^{k+1}}.
\end{equation}

\medskip

Somewhat surprisingly, $m(n)$ has an additive definition in terms of integer 
partitions. We say that a partition $(1^{a_1},2^{a_2},\dots, k^{a_k})$ of $n$ 
is 
{\em perfect}, if for every $m<n$ there is exactly one
$k$-tuple $(b_1,\dots,b_k)$, $0\le b_i\le a_i$ for all $i$, 
such that $(1^{b_1},2^{b_2},\dots, k^{b_k})$ is a partition of $m$. MacMahon 
\cite{macm91a} proved the identity 
$$
m(n)=\#\,\mbox{perfect partitions of $(n-1)$}.
$$
For example, since $m(12)=8$, 
we have 8 perfect partitions of $11$, namely 
$(1^2,3,6)$, $(1,2^2,6)$, $(1^5,6)$, $(1,2,4^2)$, $(1^3,4^2)$, $(1^2,3^3)$, 
$(1,2^5)$, 
and $(1^{11})$.

\medskip

In conclusion of the survey of previous results we should remark that from an 
enumerative point of view
it is natural to consider $m(n)$ as a function of the partition 
$\lambda=(a_1,a_2,\dots,a_k)$ of 
$\Omega(n)$, where $n=q_1^{a_1}q_2^{a_2}\dots q_k^{a_k}$ with 
$a_1\ge a_2\ge\dots\ge a_k$, rather than $n$. 
Then $m(\lambda)$ is defined as the number of ways to write 
$\lambda=v_1+v_2+\cdots+v_t$ where each 
$v_i$ is a $k$-tuple of nonnegative integers, the order of summands matters, 
and no $v_i$ is a zero vector. 
So $m(\lambda)$ is naturally understood as the number of $k$-dimensional 
compositions of $\lambda$. 
This approach was pursued by MacMahon in 
his memoirs \cite{macm91a}, \cite{macm91b}, and \cite{macm93}, see also 
\cite{macm_book}.

\medskip

The sequence
$$
(m(n))_{n\ge 1}=(1,1,1,2,1,3,1,4,2,3,1,8,1,3,3,8,1,8,1,8,3,3,1,20,2,\dots)
$$
forms entry A074206 of the database \cite{sloa}. 
Continuing the sequence a little further, we notice that $m(48)=48$ and that 
$n=48=2^4\cdot 3$ is the smallest $n>1$ such that 
$m(n)=n$. The first formula in (\ref{mforpq}) produces infinitely many $n$ 
with this property: setting $n=2^{2q-2}q$ with a prime $q>2$, we get $m(n)=n$.
We record this observation as follows:

\begin{proposition}
There exist infinitely many positive integers $n$ such that $m(n)=n$.
\end{proposition}

\noindent
This result was obtained independently also by Knopfmacher and Mays \cite{knop_mays}.

\medskip

We look at periodicity properties of the numbers $m(n)$. The recurrence 
(\ref{Moebinv}) implies easily the following result.

\begin{proposition}
The number $m(n)$ is odd if and only if $n$ is squarefree. 
\end{proposition}

It would be interesting to characterize the behavior of $m(n)$ with respect to 
other moduli besides $2$. In the next Proposition we give a partial result in this 
direction. Recall 
that an integer valued function $f(n)$ defined on the set of positive 
integers is called
{\it eventually periodic} modulo $k$ if there exist integers 
$n_0$ and $T$ such that $f(n)\equiv f(n+T)\pmod k$ for all $n>n_0$. We show that 
$m(n)$ is not eventually periodic modulo $k$ by proving a stronger result that 
$m(n)$ is not eventually constant modulo $k$ on any infinite arithmetic 
progression with coprime difference and the first term.

\begin{proposition}
The function $m(n)$ is not eventually constant modulo $k$, where $k\ge 2$, on any 
infinite arithmetic progression $n\equiv A\pmod K$, $K\ge 2$, with coprime $A$ and $K$. 
\end{proposition}

\begin{proof}
By Dirichlet's theorem, this arithmetic progression contains infinitely many prime
numbers and therefore $m(n)=1$ for infinitely many $n\equiv A\pmod K$. We select a prime
$q$ not dividing $K$ and an integer $z$ (coprime with $K$) such that $qz\equiv A\pmod K$. 
Since there are infinitely many prime numbers congruent to $z$ modulo $K$, there are 
also infinitely many $n\equiv A\pmod K$ of the form $qp$ where $p$ is a prime. Thus 
there are infinitely many $n\equiv A\pmod K$ with $m(n)=3$. Because $1\not\equiv 3\pmod k$ 
for $k>2$, we are done if $k>2$. For $k=2$, $m(n)\equiv 1\pmod 2$ for infinitely many 
$n\equiv A\pmod K$ as before. As we noted, $m(n)$ is even 
iff $n$ is not squarefree. It follows that $m(n)\equiv 0\pmod 2$ for infinitely many 
$n\equiv A\pmod K$ as well, which settles the case $k=2$.
\end{proof}

\noindent
The condition $(A,K)=1$ cannot be omitted because if $(A,K)$ is not squarefree, $m(n)$ is even 
for all $n\equiv A\pmod K$. 

\medskip

Recall now that a sequence $(f(n))_{n\ge 1}$ is holonomic 
if there exist positive integer polynomials $g_0,\dots, g_k$, not all zero, 
such that 
\begin{equation}
\label{eq:rec}
g_k(n)f(n+k)+g_{k-1}(n)f(n+k-1)+\dots+g_0(n)f(n)=0\qquad {\text{\rm 
for~all}}~n\ge 1. 
\end{equation}

\begin{proposition}
The sequence $m(n)$ is not holonomic.
\end{proposition}

\begin{proof}
Dividing (\ref{eq:rec}) by one of the (nonzero) coefficients $g_j$ with the 
largest 
degree, we obtain 
the relation 
$$
f(n+j)=\sum_{0\le i\le k,i\ne j}h_i(n)f(n+i)
$$ 
where the $h_i$'s are rational functions such that each $h_i(x)$ goes to a 
finite constant 
$c_i$ as $x\to\infty$ (we may even assume that $|c_i|\le 1$ for every $i$). 
Hence there is a constant $C>0$ (depending only on $k$ and the polynomials 
$g_i$) such that 
$$
|f(n)|\le C\max\left\{|f(n+i)|:\;-k\le i\le k,i\ne 0\right\}\ 
\mbox{ for every $n\ge k+1$.}
$$ 
We show that $(m(n))_{n\ge 1}$ violates this property.

We fix two integers $k,a\ge 1$ with the only restriction that $a$ is 
coprime to 
each of the numbers 
$1,2,\dots,k$. It is an easy consequence of the Fundamental Lemma of the 
Combinatorial Sieve
(see \cite{HR}) that there is a constant $K>0$ depending only on $k$ so that
$$
\Omega((an-k)(an-k+1)\dots(an-1)(an+1)\dots(an+k))\le K
$$
holds for infinitely many integers $n\ge 1$. For each of these $n$'s the 
$2k$ values 
$m(an+i)$, $-k\le i\le k$ and $i\ne 0$, are bounded by a constant 
(depending only on $k$) while 
the value $m(an)$ is at least $m(a)$ and can be made arbitrarily 
large by an appropriate 
selection of $a$. This contradicts the above property of holonomic sequences. 
\end{proof}

\begin{remark} The above proof can be adapted in a straightforward way to 
show that other number theoretical functions such as 
$\omega(n),~\Omega(n)$ and $\tau(n)$, where $\tau(n)$ is the number
of divisors of $n$, are not holonomic.
\end{remark}

\medskip

We present two more estimates related to the function $m(n)$.

\begin{proposition}
\label{prop:n<x}
The estimate
$$
\#\{m(n):\;n\le x\}\le \exp\(\pi\sqrt{2/\log 8}(1+o(1))(\log x)^{1/2}\)
$$
holds as $x\to\infty$.
\end{proposition}

\begin{proof}
Because $m(n)$ depends only on the partition $a_1+\dots+a_k=\Omega(n)$, 
where $n=q_1^{a_1}\dots q_k^{a_k}$
($q_1,\dots,q_k$ are distinct primes and $a_1\ge a_2\ge \dots\ge a_k>0$ are 
integers), we have that
$$
\#\{m(n):n\le x\}\le p(1)+p(2)+\cdots+p(r)\le rp(r)
$$
where $p(n)$ denotes the number of partitions of $n$ and 
$r=\max_{n\le x}\Omega(n)$. The result follows from 
$r\le\log x/\log 2$ and the classic asymptotic relation 
$p(n)\sim\exp(\pi\sqrt{2n/3})/(4n\sqrt{3})$ due to 
Hardy and Ramanujan \cite{hard_rama}.
\end{proof}

We show that the same bound on the number of distinct values of $m(n)$ holds 
when the condition $n\le x$ is replaced
with $m(n)\le x$. 

\begin{proposition}
\label{prop:f(n)<x}
The estimate
$$
\#\{m(n):\;m(n)\le x, n\ge 1\}\le 
\exp\(\pi\sqrt{2/\log 8}(1+o(1))(\log x)^{1/2}\)
$$
holds as $x\to\infty$.
\end{proposition}

\begin{proof}
As in Proposition \ref{prop:n<x} we have
$$
\#\{m(n):\;m(n)\le x, n\ge 1\}\le  p(1)+p(2)+\cdots+p(r)\le rp(r)
$$
where now $r=\max_{m(n)\le x}\Omega(n)$. By the third inequality in
Lemma~\ref{lemma:down}, $2^{r-1}=2^{\Omega(n)-1}\le m(n)\le x$ for some $n$. Thus 
$r\le 1+\log x/\log 2$ and the result 
follows as in the proof of Proposition \ref{prop:n<x} using the asymptotics of 
$p(n)$.
\end{proof}

\end{document}